\documentclass{article}    
\newcommand{\QED}{{\par\hfill$\square$\par}}

\newcommand{\Int}[2]{{\displaystyle \int_{ #1}^{ #2}}}

\newcommand{\Frac}[2]{\displaystyle{\frac{\displaystyle{#1}}{\displaystyle{#2}}}}

\newcommand{\beea}{\begin{eqnarray}}

\newcommand{\eeea}{\end{eqnarray}}

\newcommand{\ms}{\medskip\smallskip}

\newcommand{\bfe}{{\mbox{\boldmath $e$}} }

\newcommand{\0}{{\mbox{\boldmath $0$}} }

\newcommand{\Bd}{\begin{defn}\begin{rm}}
\newcommand{\EED}[1]{\end{rm}\label{defn:#1}\end{defn}}

\newcommand{\BF}{\begin{footnotesize}}

\newcommand{\EF}{\end{footnotesize}}

\newcommand{\ode}[2]{{\displaystyle \frac{\mbox{$d #1$}}{\mbox{$d #2$}}}}

\newcommand{\bi}{\begin{itemize}}

\newcommand{\ei}{\end{itemize}}

\newcommand{\ed}{\end{document}}

\newcommand{\be}{\begin{equation}}

\newcommand{\ba}{\begin{array}}

\newcommand{\ea}{\end{array}}

\newcommand{\ee}{\end{equation}}

\newcommand{\eeq}[1]{\label{eq:#1}\end{equation}}

\newcommand{\real}{{\mathbb R}}

\newcommand{\bfchi}{\mbox{\boldmath $\chi$}}

\newcommand{\bfx}{\mbox{\boldmath $x$}}

\newcommand{\bfv}{{\mbox{\boldmath $v$}} }

\newcommand{\bfu}{{\mbox{\boldmath $u$}} }

\newcommand{\bfw}{{\mbox{\boldmath $w$}} }

\newcommand{\bfA}{{\mbox{\boldmath $A$}} }

\newcommand{\bfM}{{\mbox{\boldmath $M$}} }

\newcommand{\bfB}{{\mbox{\boldmath $B$}} }

\newcommand{\bfI}{{\mbox{\boldmath $I$}} }

\newcommand{\bfN}{{\mbox{\boldmath $N$}} }

\newcommand{\cald}{{\cal D}}

\newcommand{\calp}{{\cal P}}

\newcommand{\calq}{{\cal Q}}

\newcommand{\bfK}{{\mbox{\boldmath $K$}} }

\newcommand{\bfL}{{\mbox{\boldmath $L$}} }

\newcommand{\bfg}{{\mbox{\boldmath $g$}} }

\newcommand{\bfn}{{\mbox{\boldmath $n$}} }

\newcommand{\half}{\mbox{$\frac{1}{2}$}}

\def\Bbb R{\real}

\def\tilde{\widetilde}

\def\bar{\overline}

\newcommand{\bfgamma}{\mbox{\boldmath $\gamma$}}

\newcommand{\ED}{\end{description}}

\def\tag{\renewcommand{\theequation}}

\newcommand{\Br}{\begin{rem}\begin{rm}}

\newcommand{\Er}{\end{rm}\end{remark}}

\newtheorem{lemm}{Lemma}[section]

\newtheorem{theo}{Theorem}[section]

\newtheorem{prop}{Proposition}[section]

\newtheorem{rem}{Remark}[section]

\newtheorem{defn}{Definition}[section]

\newtheorem{coro}{Corollary}[section]

\newtheorem{exe}{\footnotesize{Exercise}}[section]

\newcommand{\Be}{\begin{exe}\begin{footnotesize}\begin{rm}}

\newcommand{\EE}[1]{\end{rm}\end{footnotesize}\label{exe:#1}\end{exe}}

\newcommand{\Bt}{\begin{theo}\begin{sl}}

\newcommand{\Bp}{\begin{prop}\begin{sl}}

\newcommand{\EP}[1]{\end{sl}\label{prop:#1}\end{prop}}

\newcommand{\Et}{\end{sl}\end{theorem}}

\newcommand{\Bl}{\begin{lemm}\begin{sl}}

\newcommand{\El}{\end{sl}\end{lemma}}

\newcommand{\eqref}[1]{{\rm (\ref{eq:#1})}}

\newcommand{\Bc}{\begin{coro}\begin{sl}}

\newcommand{\Ec}{\end{sl}\end{coro}}

\newcommand{\ET}[1]{\end{sl}\label{theo:#1}\end{theo}}

\newcommand{\EL}[1]{\end{sl}\label{lemm:#1}\end{lemm}}

\newcommand{\theoref}[1]{{\rm Theorem \ref{theo:#1}}}

\newcommand{\propref}[1]{{\rm Proposition \ref{prop:#1}}}

\newcommand{\ER}[1]{\end{rm}\label{rem:#1}\end{rem}}

\newcommand{\EC}[1]{\end{sl}\label{coro:#1}\end{coro}}

\newcommand{\lemmref}[1]{{\rm Lemma \ref{lemm:#1}}}

\usepackage{graphicx}
\newcommand\smallH{
  \mathchoice
    {{\scriptstyle\mathcal{H}}}
    {{\scriptstyle\mathcal{H}}}
    {{\scriptscriptstyle\mathcal{H}}}
    {\scalebox{.7}{$\scriptscriptstyle\mathcal{O}$}}
  }
\newcommand\smallL{
  \mathchoice
    {{\scriptstyle\mathcal{L}}}
    {{\scriptstyle\mathcal{L}}}
    {{\scriptscriptstyle\mathcal{L}}}
    {\scalebox{.7}{$\scriptscriptstyle\mathcal{O}$}}
  }

\usepackage{amssymb}
\usepackage{amsfonts}

\renewcommand{\real}{{\mathbb R}}

\usepackage[mathscr]{eucal}
\usepackage{color}
\usepackage{sectsty}
\allsectionsfont{\bfseries\sffamily}
\begin{document}
\title{Stability and Long-Time Behavior of a Pendulum with an Interior Cavity Filled with a Viscous Liquid} 
\author{G.P. Galdi \thanks{Department of Mechanical Engineering and Materials Sciences, University of Pittsburgh, USA;\newline \hspace*{4mm} email: galdi@pitt.edu}\ \ \  \&\  G. Mazzone
\thanks{Department of Mathematics, Vanderbilt University, USA; email: giusy.mazzone@vanderbilt.edu}}
\date{}
\maketitle
\begin{abstract}We show  asymptotic, exponential stability of the equilibrium configuration, $\smallL$, of a hollow physical pendulum with its inner part entirely filled with a viscous liquid, corresponding to the center of mass being in the lowest position. Moreover, we prove that every weak solution with initial data possessing finite total initial energy and belonging to a ``large" open set, becomes eventually smooth and decays exponentially fast to the equilibrium $\smallL$. These results are obtained also as byproduct of a ``generalized linearization principle" that we show for evolution equations with non-empty ``slow" center manifold.  

\end{abstract}
\renewcommand{\theequation}{{1}.\arabic{equation}}
\setcounter{section}{1}
\pagenumbering{arabic}
\section*{\large Introduction} The motion of a rigid body with an interior cavity entirely filled with a viscous liquid represents one of the most important problems in the area of fluid-structure interactions. In addition to its substantial relevance in many applications, such as aerospace engineering and geophysics (see, e.g., the comprehensive monograph \cite{CAL}) it presents a number of intriguing questions  of great appeal to the applied mathematician. Thus, it is not surprising that this area of research has collected hundreds of remarkable dedicated papers and a few monographs that would be too long to list here, and for which we refer the reader to \cite[pp. x--xxxi]{CAL}.   
\par
One of the main characteristics of the evolution of such coupled systems is that the presence of the liquid can affect in a substantial way the motion of the rigid body, and may eventually produce a -sometimes unforeseen-  stabilizing effect. More precisely, after an initial ``chaotic" behavior, whose duration, $t_0$, depends on the ``size" of the initial data as well as on the relevant physical parameters involved (viscosity and density of the liquid, mass distribution of the rigid body, etc.) \cite{KL,LK},  the system reaches a steady-state, where the liquid is motionless (relative to the body), while the body executes a time-independent motion. In some cases, the latter may even reduce to an equilibrium configuration, namely, in other words, the  effect of the liquid is to bring the {\em whole} coupled system to the terminal state of rest. \par A most significant example where the latter situation occurs is a hollow physical pendulum with its inner part completely filled with a viscous liquid. The study of the motion of such a system has been carried out all along by several authors under different simplifying assumptions; we refer, among others, to \cite{kr,mr,ch,ch1}. However, it was only recently that the present authors provided a rigorous mathematical analysis of the problem on the {\em full} set of equations, without recurring to any simplification or approximation \cite{GaMa}. Their study, based on an appropriate adaptation of classical dynamic system theory, shows, among other things, that a pendulum having its interior entirely filled with a Navier-Stokes liquid reaches, eventually, an equilibrium configuration where the liquid is at relative rest and the center of mass of the coupled system body-liquid is either in its higher ($\smallH$) or lower ($\smallL$) position. This result only requires that the total initial energy of the coupled system is finite but, otherwise, of arbitrary magnitude. Moreover, the $\smallL$-configuration is stable and attainable from a ``large" class of initial data, while the $\smallH$-configuration is unstable.

At this point we wish to emphasize another not less significant phenomenon,  which mostly motivates  the writing of this article. Actually,  in both numerical and lab tests  \cite{DGMZ,WJ}, it is observed is that  after the time $t_0$ has elapsed, the system  reaches the steady-state configuration in a rather {\em abrupt} fashion. This seems to suggest that, after a sufficiently large time  that allows the liquid to be {\em almost} at rest, the whole system approaches  the terminal state at an {\em exponentially fast} rate. One of the main objectives of this paper is to show that this is indeed the case for the system constituted by a physical pendulum with an interior cavity entirely filled with a viscous liquid.        
\par
Seemingly, in order to achieve our goal, the method used in \cite{GaMa} is not particularly effective, and we have to resort to a different one. To this end, also motivated by the results obtained in \cite{Ga_im} for  questions of similar nature, in Section 1 we propose a general approach in  a class of nonlinear evolution problems with a ``slow" (local) center manifold, that is,  the spectrum, $\sigma$,  of the relevant linear (time-independent) operator, $\bfL$, is discrete with $\sigma\cap \{{\rm i}\,\real\}=\{0\}$. Actually, as directly or indirectly showed in \cite{KK,DGMZ,Ga_im,GMM}, the existence of such a manifold appears to be a {\em basic characteristic} of this kind of fluid-structure interaction problems. This is due to the fact that, for obvious physical reasons, the set of steady-state solutions does not reduce to a singleton, and may even form a continuum, either in absence or presence of a driving force. In this sense, we believe that the approach here presented might be useful also in other circumstances as, for example, those analyzed in \cite{GMM}.
\par
In addition to those mentioned above, the main assumptions we impose on $\bfL$ are, basically, that it is Fredholm of index 0, sectorial,  with $\Re[\sigma\backslash\{0\}]>0$. Moreover, null and range spaces of $\bfL$ share  only the zero element. Under these hypotheses, and some other technical ones on the nonlinear operator --compatible with the existence of multiple steady-state solutions-- we show that for sufficiently ``small" initial data the generic corresponding solution to the relevant evolution problem (see \eqref{1}) will tend to an element of the null space of $\bfL$  exponentially fast; see \theoref{1}. This result is in the spirit of  ``generalized linearization principles" like, for example, that of \cite[Theorem 2.1]{PSZ}, even though some of our assumptions and method of proof (fractional powers) are different and specifically aimed at treating the above type of fluid-structure interaction problems.  
\par     
The general theory developed in Section 1 is then  applied, in Section 2, to study the stability of $\smallL$- and $\smallH$- equilibria of the pendulum with a liquid-filled cavity. We thus show, at first, that the problem can be formulated in a suitable Hilbert space where linear and nonlinear operators satisfy all the assumptions of the theory. As a result, we prove in \theoref{2} that the $\smallL$-equilibrium is asymptotically {\em exponentially} stable, whereas the  $\smallH$-equilibrium is unstable. Successively, in the final Section 3, we address the problem of  ``abrupt" decay to (stable) equilibrium  mentioned earlier on, in the very general class of weak solution (\`a la Leray-Hopf). By combining \theoref{1} with \cite[Theorem 1.4.4]{GaMa} we thus show that for every weak solution corresponding to a ``large" set of data with finite total (kinetic and potential) energy, there exists a time, $t_0$, such that, after $t_0$, the solution becomes smooth and, together with its first time derivative and up to second spatial derivatives, must decay exponentially fast to the stable equilibrium; see \theoref{3}. 
\par  
We conclude this introductory section by some remarks about the notation used. We shall adopt standard symbols for Lebesgue, Sobolev and Bochner spaces; see, e.g., \cite{Ga}. Moreover, by the letters/symbols $c$, $C$, $c_1$, $c_2$, $C_1$, $C_2$, etc., we denote positive constants whose specific value is irrelevant, and may vary from a line to the next. If we want to emphasize the dependence on the quantity $\xi$, we shall write $c(\xi)$, etc.

\section*{\large 1. A General Approach}
Let $X$ be a (real) Banach space. We consider in $X$ the following evolution problem
\be
\ode{\bfu}t+\bfL\bfu+\bfN(\bfu)=\0\,,\ \ \bfu(0)\in X\,,
\eeq{1}
where the involved operators satisfy certain appropriate conditions that we are about to state. To this end, let $\bfA:X\mapsto X$ be a linear, sectorial operator with compact inverse and $\Re[\sigma(\bfA)]>0$. For $\alpha\in [0,1]$, set 
$$
X_{\alpha}=\big\{\bfu\in X:\ \|\bfu\|_\alpha:=\|\bfA^\alpha\bfu\|<\infty\big\}\,;\ \ X_0\equiv X\,,\ \ \|\bfu\|_0\equiv \|\bfu\|\,.
$$  
It is well known that $X_\alpha$ is a Banach space that, in addition, is compactly embedded in $X$ for $\alpha>0$, e.g., \cite[Theorem 1.4.8]{Henry}. Next, let $\bfB:X\mapsto X$ be a linear operator with $D(\bfB)\supset D(\bfA)$, and such that
\be
\|\bfB\bfu\|\le c_1\,\|\bfu\|_\alpha\,,\ \ \alpha\in [0,1)\,.
\eeq{0}
We now define the operator $\bfL$ in \eqref{1} by setting
\be 
\bfL=\bfA+\bfB\,,
\eeq{2}
with $D(\bfL)\equiv D(\bfA)$.
Since 
\be
\|\bfu\|_\alpha\le c\,\|\bfu\|^{1-\alpha}\|\bfA\bfu\|^\alpha\,,
\eeq{2i}
and $\bfA^{-1}$ is compact, it follows  that $\bfB$ is $\bfA-$compact, so that  $\bfL$ is an unbounded Fredholm operator of index 0 \cite[Theorem 4.3]{GG}. Also, 
from \cite[Theorem 1.3.2]{Henry} it follows, in particular, that $\bfL$ is sectorial. Finally, observing that by \eqref{2i},  for any $\varepsilon>0$,
$$ 
\|\bfB\bfu\|\le c(\varepsilon)\|\bfu\|+\varepsilon\,\|\bfA\bfu\|
$$
and that, by the properties of $\bfA$,
$$
\|(\lambda-\bfA)^{-1}\|\le c_2\lambda^{-1}\,,\ \ \mbox{all $\lambda>0$}
$$
we deduce \cite[Theorem 3.17 at p. 214]{Kato} that $\bfL$ has a compact resolvent and, therefore, a discrete spectrum. 
\par
We shall now
make the  following further assumptions on the operator $\bfL$:
\tag{H1}
\be
{\rm dim}\,{\sf N}[\bfL]=m\ge1\,,\vspace*{-2mm}
\label{eq:H1}\ee
\tag{H2}
\be
{\sf N}[\bfL]\cap {\sf R}[\bfL]=\{\0\}\,,
\label{eq:H2}\ee and
\tag{H3}
\be
\sigma(\bfL)\cap\{{\rm i}\,\real\}=\{0\}\,.
\label{eq:H3}
\ee\setcounter{equation}{3}
\renewcommand{\theequation}{{1}.\arabic{equation}}\smallskip\par
A first important consequence of some of the stated properties of $\bfL$ is derived next.
\Bl The space $X$ admits the following decomposition
\be
X={\sf N}[\bfL]\oplus{\sf R}[\bfL]\,.
\eeq{3}
Moreover, denoting by $\calq$ and $\calp$ the spectral projections according to the spectral sets
$$\sigma_0(\bfL):=\{0\}\,,\ \
\sigma_1(\bfL):= \sigma(\bfL)\backslash\sigma_0(\bfL)\,,
$$
we have
\be
{\sf N}[\bfL]=\calq(X)\,,\ \ {\sf R}[\bfL]:=\calp(X)\,.
\eeq{4}
Finally, \eqref{3} completely reduces $\bfL$ into $\bfL=\bfL_0\oplus\bfL_1$
with
\be
\bfL_0:=\calq\bfL=\bfL\calq\,,\ \ \bfL_1:=\calp\bfL=\bfL\calp\,, 
\eeq{5}
and $\sigma(\bfL_0)\equiv\sigma_0(\bfL)$, $\sigma(\bfL_1)\equiv\sigma_1(\bfL)$.
\EL{1}
{\em Proof.} Since $\bfL$ is Fredholm of index 0, from (\ref{eq:H1}) we deduce ${\rm codim}({\sf R}[\bfL])=m$. Thus, there exists at least one $S\textcolor{black}{\subset} X$ such that $X=S\oplus{\sf R}[\bfL]$, with $S\cap {\sf R}[\bfL]=\{\0\}$.  However,  ${\rm dim}\,(S)={\rm dim}\,({\sf N}[\bfL])=m$ and (\ref{eq:H2}) holds, so that we may take $S={\sf N}[\bfL]$, which proves \eqref{3}. The remaining properties stated in the lemma are then a consequence of \eqref{3} and classical results on spectral theory (e.g., \cite[Proposition A.2.2]{Lun}, \cite[Theorems 5.7-A,B]{Tay})
\QED\par
We now turn to the assumptions needed on the operator $\bfN$. We begin to require the following Lipschitz-like condition
\tag{H4}
\be
\|\bfN(\bfu_1)-\bfN(\bfu_2)\|\le c\,\|\bfu_1-\bfu_2\|_\alpha\,,\ \ \mbox{for all $\bfu_1,\bfu_2$ in a neighborhood of $\0\in X$}\,.
\label{eq:H4}
\ee
Furthermore, we observe that, by \eqref{3},  every $\bfu\in X$ can be written as
$$
\bfu=\bfu^{(0)}+\bfu^{(1)}\,,\ \ \bfu^{(0)}\in {\sf N}[\bfL]\,,\ \bfu^{(1)}\in {\sf R}[\bfL]\,.
$$
Thus, setting
$$
\bfM(\bfu^{(0)},\bfu^{(1)}):=\bfN(\bfu^{(0)}+\bfu^{(1)})
$$
we suppose
\tag{H5}
\be
\|\bfM(\bfu^{(0)},\bfu^{(1)})\|\le c\,\left[(\|\bfu^{(0)}\|+\|\bfu^{(1)}\|^{\kappa_1})\|\bfu^{(1)}\|^{\kappa_2}+\|\bfu^{(1)}\|_\alpha^{\kappa_3}\right]\,,\ \ \kappa_1,\kappa_2\ge 1,\ \kappa_3>1\,.
\label{eq:H5}
\ee
\renewcommand{\theequation}{1.\arabic{equation}}\setcounter{equation}{6}
\Br
Because of \eqref{H5}, it follows that ${\sf N}[\bfL]$ is contained in the set of equilibria (steady-state solutions) to \eqref{1}.  
\ER{1}

We are now in a position to prove the following stability result.
\Bt Suppose the operators $\bfL$, defined in \eqref{2}, and $\bfN$ satisfy  hypotheses {\rm (\ref{eq:H1})--(\ref{eq:H5})}. Then, if $\Re[\sigma(\bfL_1)]>0$, we can find $\rho_0>0$ such that if
$$
\|\bfu(0)\|_\alpha<\rho_0\,,
$$
there is a unique corresponding solution $\bfu=\bfu(t)$ to \eqref{1} for all $t>0$, satisfying
\be
\bfu\in C([0, T];X_\alpha)\cap C((0,T]; X_1)\cap C^1((0,T];X)\,,\ \ \mbox{for all $T>0$}\,.
\eeq{6}  
Moreover, the solution $\bfu=\0$ to \eqref{1} is exponentially stable in $X_\alpha$, namely, the following properties hold. 
\begin{itemize}
\item[{\rm (a)}] For any $\varepsilon>0$ there is $\delta>0$ such that
$$
\|\bfu(0)\|_\alpha<\delta\ \ \Longrightarrow\ \ \sup_{t\ge 0}\|\bfu(t)\|_\alpha <\varepsilon\,;
$$ 
\item[{\rm (b)}] There are $\eta,c,\kappa>0$ such that 
$$  
\|\bfu(0)\|_\alpha<\eta\ \ \Longrightarrow\ \ \|\bfu(t)-\bar{\bfu}\|_\alpha\le c\,\|\bfu^{(1)}(0)\|_\alpha\,{\rm 
e}^{-\kappa\,t}\,,\ \mbox{all $t>0$}\,,
$$
for some $\bar{\bfu}\in {\sf N}[\bfL]$.
\end{itemize}
Finally, if $\Re[\sigma(\bfL_1)]\cap (-\infty,0)\neq\emptyset$ then the solution $\bfu=\0$ to \eqref{1} is unstable in $X_\alpha$, namely, the property provided in {\rm (a)} does not hold.
\ET{1}
{\em Proof.} Under the stated assumptions on $\bfA$, $\bfB$ and (\ref{eq:H4}), the claimed instability property follows from \cite[Theorem 5.1.3]{Henry}. Likewise, the existence of a unique solution  $\bfu$ to \eqref{1} in some time interval $(0,t_\star)$ satisfying \eqref{6} for each $T\in (0,t_\star)$ is guaranteed, under the above assumptions, by classical results on semilinear evolution equations (e.g., \cite[p. 196--198]{Pazy}).
Furthermore, this solution can be extended to provide a solution beyond any time $\tau\in [T,t_\star)$ if $\|\bfu(\tau)\|_\alpha<\infty$, whereas, if $t_\star<\infty$, it will fail if and only if $\lim_{t\to t^*}\|\bfu(t)\|_\alpha=\infty$. 
We shall next show that, in fact,  only the former situation occurs, if the size of the initial data is suitably restricted.
Applying  $\calq$ and $\calp$  on both sides of \eqref{1} 
 and taking into account \eqref{5} we  show 
\be\ba{rl}\medskip
\ode{\bfu^{(1)}}t+\bfL_1\bfu^{(1)}=&\!\!\!-\mathcal P\bfM(\bfu^{(0)}, \bfu^{(1)})\\
\ode{\bfu^{(0)}}t =&\!\!\!-\mathcal Q\bfM(\bfu^{(0)}, \bfu^{(1)})\,.
\ea
\eeq{1.20}
Since the operator $\bfL$, being sectorial, is the generator of an analytic semigroup in $X$,  so is $\bfL_1$ in $X^{(1)}\equiv{\sf R}[\bfL]$.  Thus, for all $t\in[0,t_\star)$ from \eqref{1.20}$_1$ we have
\be
\bfu^{(1)}(t)={\rm e}^{-\bfL_1 t}\bfu_0^{(1)}-\int_0^t{\rm e}^{-\bfL_1 (t-s)}[\mathcal P\bfM(\bfu^{(0)}(s), \bfu^{(1)}(s))]ds\,.
\eeq{1.22}
Also, by assumption and the spectral property of $\bfL$,  there is $\gamma>0$ such that 
\be
\Re[\sigma(\bfL_1)]>\gamma>0\,,  
\eeq{1.17}
which implies that the fractional powers $\bfL_1^\alpha$, $\alpha\in (0,1)$, are well defined in $X^{(1)}$. Thus, setting
$$
\bfw:={\rm e}^{bt}\bfL_1^\alpha\bfu^{(1)}\,,\ \ 0<b<\gamma\,,
$$
from \eqref{1.22} we get
\be
\bfw(t)={\rm e}^{bt}{\rm e}^{-\bfL_1 t}\bfL_1^\alpha\bfu_0^{(1)}-\int_0^t{\rm e}^{bt}\bfL_1^{\alpha} {\rm e}^{-\bfL_1 (t-s)}[\mathcal P\bfM(\bfu^{(0)}(s), {\rm e}^{-bs}\bfL_1^{-\alpha}\bfw(s)]ds\,.
\eeq{1.23}
In view of the stated properties of $\bfL_1$ it results (e.g. \cite[Theorem 1.4.2]{Henry}), 
$$
\|\bfL_1^{-\alpha}\bfw\|\le c_1\,\|\bfw\|\,,\ \  \mbox{for all $\bfw\in X^{(1)}$}\,,
$$
and so, by (\ref{eq:H5}) and the latter, we derive
\be
\|\mathcal P\bfM(\bfu^{(0)}(s),{\rm e}^{-bt}\bfL_1^{-\alpha}\bfw(s))\|\le c_2
\,\big[(\|\bfw\|+\|\bfu^{(0)}\|)^{\kappa_1}\|\bfw\|^{\kappa_2}+\|\bfw\|^{\kappa_3}\big]\,.
\eeq{1.29}
Next, we recall that in $X^{(1)}$ it is 
\be\|\bfL_1^{\alpha}{\rm e}^{-\bfL_1t}\|\le t^{-\alpha}{\rm e}^{-\gamma t}\,, 
\eeq{df0} 
and observe that from \eqref{0} and 
 \cite[Theorem 1.4.6]{Henry} 
\be
\|\bfL_1^\alpha\bfu_0^{(1)}\|\le c_3\,\|\bfu_0^{(1)}\|_\alpha\,.
\eeq{df1}
Thus collecting  \eqref{1.23}--\eqref{df1} we deduce
\be\ba{ll}\medskip
\|\bfw(t)\|\le {\rm e}^{-(\gamma-b)t}\|\bfu_0^{(1)}\|_\alpha\\
\hspace*{2cm}+c_4\Int0t\Frac{{\rm e}^{-(\gamma-b)(t-s)}}{(t-s)^\alpha}[(\|\bfw(s)\|+\|\bfu^{(0)}(s)\|)^{\kappa_1}\|\bfw(s)\|^{\kappa_2}+\|\bfw(s)\|^{\kappa_3}]\,.\ea
\eeq{1.30}
From the local existence theory considered earlier on, we know that for any given $\rho>0$ there exists an interval of time $[0,\tau]$, $\tau<t_\star$, such that 
\be
\sup_{t\in [0,\tau)}\left(\|\bfw(t)\|+\|\bfu^{(0)}(t)\|\right)\le \rho\,,\ \tau<t_\star\,
\eeq{71}
provided 
$\|\bfu(0)\|_\alpha<\eta$, for some $\eta>0$. 
Our first objective is to show that  $\eta$ and $\rho$ can be chosen sufficiently small so that \eqref{71} holds also with $\tau=t_\star$, thus implying, in particular, that the solution $\bfu=\bfu(t)$ to \eqref{1} exists for all times $t>0$. In fact, suppose, 
by contradiction, that there is $\tau_0<t_\star$ such that
\be
\|\bfw(t)\|+\|\bfu^{(0)}(t)\|<\rho\,,\ \ t\in [0,\tau_0)\ \ \mbox{and}\ \ \|\bfw(\tau_0)\|+\|\bfu^{(0)}(\tau_0)\|=\rho\,.
\eeq{72}
Noticing that
$$
\int_0^t\frac{{\rm e}^{-(\gamma-b)(t-s)}}{(t-s)^\alpha}\,ds\le \int_0^\infty\frac{{\rm e}^{-(\gamma-b)t}}{t^\alpha}\,dt<\infty\,,
$$
from \eqref{1.30} and \eqref{72} we get
\be
\|\bfw(t)\|\equiv {\rm e}^{bt}\|\bfu^{(1)}(t)\|_\alpha\le \eta + \epsilon(\rho)\, \rho\,,\ \mbox{for all $t\in [0,\tau_0]$},
\eeq{1.31}
where $\epsilon(\rho)$ represents, here and in what follows, a generic smooth, positive function such that $\epsilon(\rho)\to 0$ as $\rho\to 0$.
On the other hand, \eqref{1.20}$_2$  with the help of (\ref{eq:H5}), furnishes
$$
\|\bfu^{(0)}(t)\|\le \|\bfu^{(0)}(0)\|+c_5\int_0^t \,\big[(\|\bfu^{(1)}(s)\|+\|\bfu^{(0)}(s)\|)^{\gamma_1}\|\bfu^{(1)}(s)\|^{\gamma_2}+\|\bfu^{(1)}(s)\|_\alpha^{\gamma_3}\big]\,ds\,.
$$
Thus, if we restrict ourselves to $t\in [0,\tau_0]$,  and use   \eqref{1.31} and \eqref{72}, the preceding inequality provides
\be\ba{ll}\medskip
\|\bfu^{(0)}(t)\|\le &\!\!\!\|\bfu^{(0)}(0)\|+\epsilon(\rho)\Int0t \,\|\bfu^{(1)}(s)\|_\alpha\,ds\le \|\bfu^{(0)}(0)\|+\epsilon(\rho)\Int0t \,{\rm e}^{-b\,s}\|\bfw(s)\|\,ds \\&\!\!\!\!\le 2\eta +\epsilon(\rho)\,\rho\,,\ \ t\in [0,\tau_0]\,. \ea
\eeq{1.32} 
Combining \eqref{1.31} and \eqref{1.32}, and choosing $2\eta/\rho+ \epsilon(\rho)<1/4$ we  conclude in particular
$$
\|\bfw(\tau_0)\|+\|\bfu^{(0)}(\tau_0)\|\le \rho/2
$$
contradicting \eqref{72}. As a result, by what we observed early on, we may take $t_*=\infty$ in \eqref{71} and conclude as well
\be
\sup_{t\in [0,\infty)}\left(\|\bfw(t)\|+\|\bfu^{(0)}(t)\|\right)\le \rho\,,
\eeq{711}
proving, as a byproduct, the desired global existence property.
Now, from \eqref{711}  and \eqref{1.30}, we easily deduce, for $\rho$ small enough, 
\be
\|\bfw(t)\|\le c_6\,\|\bfu^{(1)}(0)\|_\alpha\,,\ \ \mbox{all $t>0$}\,,
\eeq{sfc}
namely,
\be
\|\bfu^{(1)}(t)\|_\alpha\le c_7 \,{\rm e}^{-b t} \|\bfu^{(1)}(0)\|_\alpha\,,\ \ \mbox{all $t>0$}\,.  
\eeq{1.35}
Also, employing \eqref{sfc} into \eqref{1.32}, we infer
\be
\|\bfu^{(0)}(t)\|\le c_8\,\|\bfu(0)\|_\alpha\,,\ \ \mbox{all $t>0$}\,.
\eeq{sfc1}
Therefore, from \eqref{1.35} and \eqref{sfc1} we recover the stability property stated in (a).
Moreover, integrating  \eqref{1.20}$_2$ between arbitrary $t_1,t_2>0$ using \eqref{1.35} and reasoning in a way similar to what we did to obtain \eqref{1.32} we get
\be
\|\bfu^{(0)}(t_1)-\bfu^{(0)}(t_2)\|\le c_{9}\int_{t_1}^{t_2}\|\bfu^{(1)}(s)\|_\alpha ds\le c_{10}\,\|\bfu^{(1)}(0)\|_\alpha\int_{t_1}^{t_2}{\rm e}^{-b s}ds\,, 
\eeq{1.36}
from which we deduce that there exists $\bar{\bfu}\in {\sf N}[\bfL]$ such that
$$
\lim_{t\to\infty}\|\bfu^{(0)}(t)-\bar{\bfu}\|=0\,.
$$
Employing this information into \eqref{1.36} in the limit $t_2\to\infty$, and with $t_1=t$ we show
$$
\|\bfu^{(0)}(t)-\bar{\bfu}\|_\alpha\le c_{11}\,\|\bfu^{(0)}(t)-\bar{\bfu}\|\le c_{12}\|\bfu^{(1)}(0)\|_\alpha\,{\rm e}^{-b t},
$$
which, once combined with \eqref{1.35}, proves the exponential rate of decay stated in (b).   The proof of the theorem is thus concluded.
\QED
\renewcommand{\theequation}{2.\arabic{equation}}\setcounter{equation}{0}\setcounter{section}{2}
\section*{\large 2. Asymptotic Stability of the Equilibrium Configurations of a Pendulum with a Liquid-Filled Cavity}
Objective of this and the next section is to apply the general theory developed in the previous one, also combined with the findings of \cite{GaMa},  to the study of the stability of the equilibria, and long-time behavior of the generic motion of a physical pendulum with an interior cavity entirely filled with a {\em viscous} liquid. \par
More specifically, let $\mathscr S$ be the coupled system  constituted by a rigid body, $\mathscr B$, with an interior cavity, $\mathscr C$ (assumed to be a domain of $\real^3$ of class $C^2$), entirely filled with a Navier-Stokes liquid, $\mathscr L$. Suppose that $\mathscr B$ is constrained to move (without friction) around a horizontal axis ${\sf a}$ in such a way that during all possible motions of $\mathscr S$ its center of mass $G$ belongs to a fixed vertical plane  orthogonal to ${\sf a}$, so that the distance from $G$ to its orthogonal projection, $O$, on ${\sf a}$ is kept constant. 

Denoting by $\mathscr F\equiv\{O,\bfe_1,\bfe_2,\bfe_3\}$  a frame attached to $\mathscr B$, with the origin at $O$, $\bfe_1\equiv\, \stackrel{\rightarrow}{O\,G}/|\stackrel{\rightarrow}{O\,G}|$, and $\bfe_3$ directed along   ${\sf a}$, we then  have that the motion of $\mathscr S$ in $\mathscr F$ is governed by the following set of equations \cite{mr}
\be
\ba{cc}\medskip\left.\ba{rl}\ms
\rho\bigl({\bfv_t}+\dot{\omega}\bfe_3\times\bfx+\bfv\cdot\nabla\bfv
+2\omega\,\bfe_3\times\bfv\bigr)&\!\!\!=\mu\Delta\bfv-\nabla p\\
\nabla\cdot\bfv&\!\!\!=0\ea\right\}\ \ \mbox{in $\mathscr C\times \real_+$}\\ \ms
\bfv(x,t)|_{\partial \mathscr C}=\0\\
{\sf C}(\dot{\omega}-\dot{a})=\beta^2\chi_2\,,\ \ \ 
\dot{\bfchi}+\omega\,\bfe_3\times\bfchi=\0\,,
\ea
\eeq{2.1}
Here, $\bfv$ and $p$ are {\em relative} {velocity} and (modified) pressure fields of $\mathscr L$, respectively, while $\rho$ and $\mu$ are its density and shear viscosity coefficient. Also, $\omega\,\bfe_3$ is the angular velocity of $\mathscr B$ and $\bfchi=(\chi_1\equiv \cos\varphi,\chi_2\equiv-\sin\varphi,0)$ where $\varphi$ is the angle between $\bfe_1$ and the gravity $\bfg$. Furthermore, ${\sf C}$ is the moment of inertia of $\mathscr S$ with respect to ${\sf a}$,
\be a:=-\frac{\rho}{\sf C}\,\bfe_3\cdot\int_{\mathscr C}\bfx\times\bfv\,,
\eeq{2.2}
and $$\beta^2=M\,g\,|\stackrel{\rightarrow}{O\,G}|\,,$$
with $M$ mass of $\mathscr S$.
\par
It is not difficult to show (formally) that \eqref{2.1} has  only two  steady-state solutions given by
\be
{\sf s}_0^{\pm}:=(\bfv\equiv\nabla p\equiv\0, \omega\equiv0,\bfchi=\pm\bfe_1)\,,
\eeq{2.3} 
and representing the equilibrium configurations where $\mathscr S$ is at rest with $G$ in its lowest (${\sf s}_0^+$) or highest (${\sf s}_0^-$) position. 

We shall next employ \theoref{1} to investigate the stability property of the above equilibrium configurations. Successively, combining that theorem with some of the results established in \cite{GaMa}, we will characterize the asymptotic behavior of the solutions to \eqref{2.1}  in a very general class, and for a ``large"  set of initial data.  
\par
To accomplish all the above,  we begin to  observe that the ``perturbed motion" around ${\sf s}_0^\pm$ can be written as
\be
(\bfv, p,\omega, \bfchi:=\bfgamma\pm\bfe_1)\,,\ \ |\bfgamma\pm\bfe_1|=1\,,
\eeq{fan}
where, by  \eqref{2.1},  the ``perturbation" $(\bfv, p,\omega, \bfgamma)$ satisfies the following equations  
\be
\ba{cc}\medskip\left.\ba{rl}\ms
\rho\bigl({\bfv_t}+\dot{\omega}\bfe_3\times\bfx\bigr)-\mu\Delta\bfv+\nabla p
&\!\!\!=
-\rho\bigl(2\omega\,\bfe_3\times\bfv+\bfv\cdot\nabla\bfv\bigr)\\
\nabla\cdot\bfv&\!\!\!=0\ea\right\}\ \ \mbox{in $\mathscr C\times \real_+$}\\ \ms
\bfv(x,t)|_{\partial \mathscr C}=\0\\
{\sf C}(\dot{\omega}-\dot{a})=\beta^2\gamma_2\,,\ \ \ 
\dot{\bfgamma}+\omega\,\bfe_3\times\bfgamma_0=-\omega\,\bfe_3\times\bfgamma\,,
\ea
\eeq{2.4}
with
\be
\bfgamma_0:=\xi\,\bfe_1\,,\ \ \xi=\pm 1\,.
\eeq{2.4_1}
Our subsequent step is to write \eqref{2.4} as an evolution problem of  the type \eqref{1} with $X$ appropriate Hilbert space $H$. We thus define
$$
L^2_\sigma(\mathscr C):=\left\{\bfv\in L^2(\mathscr C):\ \nabla\cdot\bfv=0\ \mbox{in $\mathscr C$}\,,\ \bfv\cdot\bfn|_{\partial\mathscr C}=0\right\}\,,
$$
and let
$$
H:=\bigl\{\bfu:=(\bfv,\omega,\bfgamma)^\top: \bfu\in L^2_\sigma(\mathscr C)\oplus\real\oplus\real^2\bigr\}\,,
$$
endowed with the scalar product
$$
\langle \bfu_1,\bfu_2\rangle:=\int_{\mathscr C}\bfv_1\cdot\bfv_2\,{\rm d}\mathscr C+\omega_1\,\omega_2+\bfgamma_1\cdot\bfgamma_2\,,
$$
and corresponding norm
$$
\|\bfu\|:=\langle\bfu,\bfu\rangle^{\frac12}\,.
$$
Moreover,  we introduce the following operators
\be\ba{ll}\medskip
\bfI:\bfu\in H\mapsto \bfI\bfu:=\bigl(\rho\,\bfv+{\rm P}[\rho\,\omega\bfe_3\times\bfx], {\sf C}(\omega-a),\bfgamma)^\top\in H\\ \medskip
\tilde{\bfA}:\bfu\in D(\tilde{\bfA}):=W^{2,2}(\mathscr C)\cap\cald_0^{1,2}(\mathscr C)\oplus\real\oplus\real^2\subset H\mapsto\tilde{\bfA}\bfu:=\bigl(-\mu\,{\rm P}\Delta\bfu,\omega,\bfgamma\bigr)^\top\in H\\ \medskip
\tilde{\bfB}:\bfu\in H\mapsto\tilde{\bfB}\bfu:=\bigl(\0,-\beta^2\gamma_2-\omega, \omega\bfe_3\times\bfgamma_0-\bfgamma\bigr)^\top\in H\\
\tilde{\bfN}:\bfu\in D(\tilde{\bfA})\subset H\mapsto\tilde{\bfN}(\bfu):= \bigl(-\rho{\rm P}[2\omega\bfe_3\times\bfv+\bfv\cdot\nabla\bfv],0,-\omega\bfe_3\times\bfgamma\bigr)^\top\in H\,,
\ea
\eeq{2.5}
where ${\rm P}:L^2(\mathscr C)\mapsto L^2_\sigma(\mathscr C)$ is the Helmholtz projector. It is readily shown that the operator $\bfI$ is boundedly invertible. In fact, we have
\be
\langle\bfI\bfu,\bfu\rangle=\int_{\mathscr C}\rho\,|\bfv|^2+2\rho\,\omega\int_{\mathscr C}\bfe_3\times\bfx\cdot\bfv+{\sf C}\omega^2+|\bfgamma|^2\,,
\eeq{2.6}
and so, observing that ${\sf C}={\sf C}_{\mathscr B}+{\sf C}_{\mathscr L}$ with ${\sf C}_{\mathscr B}$ and ${\sf C}_{\mathscr L}$ moments of inertia with respect to ${\sf a}$ of $\mathscr B$ and $\mathscr C$, respectively, and that
\be
{\sf C}_{\mathscr L}=\int_{\mathscr C}\rho\,(\bfe_3\times\bfx)^2\,,
\eeq{df}
from \eqref{2.6} we deduce
$$
\langle\bfI\bfu,\bfu\rangle=\int_{\mathscr C}\rho\,|\bfv+\omega\bfe_3\times\bfx|^2+C_{\mathscr B}\omega^2+|\bfgamma|^2\,.
$$
From this relation it follows that ${\sf N}[\bfI]=\{\0\}$. On the other hand, by its very definition we have $\bfI=\bfI_1+\bfK$ where $\bfI_1$ is an isomorphism of $H$ onto itself and $\bfK$ is a compact (finite-dimensional) operator. We thus conclude the stated invertibility property of $\bfI$. From \eqref{2.5}$_1$, it also follows at once that the operator $\bfI$ is  symmetric, namely,
$$
\langle\bfI\bfu_1,\bfu_2\rangle=\langle\bfu_1,\bfI\bfu_2\rangle\,,\ \ \mbox{for all $\bfu_1,\bfu_2\in H$}\,,
$$
and, therefore, so is its inverse.
Thus, defining
\be
\bfA:=\bfI^{-1}\tilde{\bfA}\,,\ \ \bfB:=\bfI^{-1}\tilde{\bfB}\,, \ \ \bfL:=\bfA+{\bfB}\,,\ \ \bfN:=-\bfI^{-1}\tilde{\bfN}\,,
\eeq{jm}
we see that \eqref{2.4} can be written as the following evolution equation in the space $H$
\be
\ode{\bfu}t+\bfL\bfu+\bfN(\bfu)=\0\,,\ \ \ \bfu(0)\in H\,.
\eeq{2.7}
We shall next show that the operators $\bfL$ and $\bfN$ defined above satisfy all the assumptions (\ref{eq:H1})--(\ref{eq:H5}) stated in the previous section. In this regard we commence to notice that, by the properties of the Stokes operator  
$$
\bfA_0:=-\mu\,{\rm P}\Delta
$$ 
with domain $D(\bfA_0):=W^{2,2}(\mathscr C)\cap \cald_0^{1,2}(\mathscr C)$ and range  $L^2_\sigma(\mathscr C)$,
it follows that $\tilde{\bfA}$ has a compact inverse and, therefore, a purely discrete spectrum which, in addition, lies on the positive real axis. Since $\bfI^{-1}$ is symmetric (and bounded), the operator ${\bfA}$ enjoys the same stated properties as $\tilde{\bfA}$. Furthermore,  $\bfB$ is bounded and therefore satisfies \eqref{0} with $\alpha=0$, and we conclude that $\bfL$ is Fredholm of index 0. We shall now prove the validity of the other assumptions.  
\subsubsection*{Proof of (\ref{eq:H1})} The equation $\bfL\bfu=\0$ in $H$ is equivalent to the following system of equations:
$$\ba{ll}\medskip
-\mu\,{\rm P}\Delta\bfv=-\nabla p\,,\ \nabla\cdot\bfv=0\,,\ \ \bfv|_{\partial\mathscr C}=\0\\  
\beta^2\gamma_2=0\,,\ \ \omega\,\bfe_3\times\bfgamma_0=\0\,,
\ea
$$
whose solutions are of the form
\be
\bfu^{(0)}=(\0,0,\sigma\bfe_1)^\top\,,\ \ \sigma\in\real\,.
\eeq{2.8}
Therefore,
$$
{\rm dim}\,{\sf N}[\bfL]=1\,.
$$\QED
\subsubsection*{Proof of (\ref{eq:H2})} The equation $\bfL\bfu=\bfu^{(0)}$, with $\bfu^{(0)}$ given in \eqref{2.8} for some $\sigma\in\real$, is equivalent to the system of equations
$$\ba{ll}\medskip
-\mu\,{\rm P}\Delta\bfv=-\nabla p\,,\ \nabla\cdot\bfv=0\,,\ \ \bfv|_{\partial\mathscr C}=\0\\  
\beta^2\gamma_2=0\,,\ \ \omega\,\bfe_3\times\bfgamma_0=\sigma\bfe_1\,.
\ea
$$
However, with the help of \eqref{2.4_1} we infer at once $\sigma=0$, which proves the claim.\QED
\subsubsection*{Proof of (\ref{eq:H3})} The equation
$$
\bfL\bfu={\rm i}\,\zeta\,\bfu\ \ \zeta\in\real\backslash \{0\}\,,
$$
is equivalent to the following ones
\be
\ba{ll}\medskip
-\mu\,{\rm P}\Delta\bfv=-\nabla p+{\rm i}\,\zeta\,\rho\,(\bfv+\omega\,\bfe_3\times\bfx)\,,\ \nabla\cdot\bfv=0\,,\ \ \bfv|_{\partial\mathscr C}=\0\\  
\beta^2\gamma_2={\rm i}\,\zeta\,{\sf C}\,(\omega-a)\,,\ \ \omega\,\bfe_3\times\bfgamma_0={\rm i}\,\zeta\,\bfgamma\,.
\ea
\eeq{2.10}
After dot-multiplying \eqref{2.10}$_1$ by $\bar{\bfv}$ ($``\overline{\,\,\,}"=\textrm{complex conjugate}$), \eqref{2.10}$_5$ by $(\beta^2/\xi)\bar{\bfgamma}$ (see \eqref{2.4_1}),  \eqref{2.10}$_4$ by $\bar{\omega}$, and integrating by parts over $\mathscr C$ as necessary, we get 
$$\ba{rl}\medskip
\mu\|\nabla\bfv\|_2^2+{\rm i}\,\zeta\,\rho\,\|\bfv\|_2^2&\!\!\!={\rm i}\,\zeta\,{\sf C}\,a\,\bar\omega\\ \medskip
-{\rm i}\,\zeta\,{\sf C}\,|\omega|^2&\!\!\!=-\beta^2\gamma_2\bar{\omega}-{\rm i}\,\zeta\,{\sf C}\, a\,\bar{\omega}\\
-{\rm i}\,\zeta\,\Frac{\beta^2}{\xi}|\bfgamma|^2&\!\!\!=\beta^2\gamma_2\bar{\omega}\,.
\ea
$$
Adding side-by-side the  displayed equations, we find
$$
\mu\|\nabla\bfv\|_2^2+{\rm i}\,\zeta\,\bigl[\|\bfv\|_2^2-{\sf C}\,|\omega|^2-\frac{\beta^2}{\xi}|\bfgamma|^2\bigr]=0\,,
$$
which implies $\bfv\equiv\0$. As a result, by using the curl operator $\nabla\times$ on both sides of \eqref{2.10}$_1$ with $\bfv=\0$, we show $\omega=0$, which, in turn, once replaced in \eqref{2.10}$_4$ furnishes $\gamma_2=0$. Thus, finally, from \eqref{2.10}$_5$ we conclude $\bfgamma=\0$, and the proof is completed. \QED
\subsubsection*{Proof of (\ref{eq:H4})} From \eqref{2.5} and \eqref{jm}, we easily show
\be\ba{rl}\medskip
\|\bfN(\bfu_1)-\bfN(\bfu_2)\|\le c\,\bigl[|\omega_1|(\|\bfv_1-\bfv_2\|_2+&\!\!\!\!|\bfgamma_1-\bfgamma_2|)+(\|\bfv_2\|_2+|\bfgamma_2|)|\omega_1-\omega_2|\\
&+\|{\rm P}(\bfv_1\cdot\nabla\bfv_1-\bfv_2\cdot\nabla\bfv_2)\|_2\bigr]
\ea
\eeq{2.12}
Next,  we observe that the fractional powers of $\tilde{\bfA}$ are given by
\be
\tilde{\bfA}^\alpha\bfu=\bigl(\bfA_0^\alpha\bfv,\omega,\bfgamma)\,,\ \ \alpha\in(0,1)\,,
\eeq{JH}
and that, being $\bfI^{-1}$ bounded, 
by Heinz inequality 
\be
 c_1\,\|\bfA^\alpha\bfu\|\le\|{\tilde{\bfA}}^\alpha\bfu\|\le c_2\,\|\bfA^\alpha\bfu\|\,.
\eeq{2.13}
Furthermore, by a classical result  \cite[Lemma 3]{KaFu}, for any $\alpha\in [\frac34,1]$ it is
\be
\|{\rm P}(\bfv_1\cdot\nabla\bfv_1-\bfv_2\cdot\nabla\bfv_2)\|_2\le c_1\,(\| 
\bfA_0^{\alpha}\bfv_1\|_2+\| 
\bfA_0^{\alpha}\bfv_2\|_2)\,\|\bfA_0^\alpha(\bfv_1-\bfv_2)\|_2.
\eeq{2.14}
The claimed property about the validity of \eqref{H4} is then a consequence of \eqref{2.12}--\eqref{2.14}.\QED

\subsubsection*{Proof of (\ref{eq:H5})} In view of \eqref{2.8} and recalling the decomposition \eqref{3}, we have
$$
\bfu=\bfu^{(0)}+\bfu^{(1)}\equiv (\0,0,\gamma_1\bfe_1)^\top+(\bfv,\omega,\gamma_2\bfe_2)^\top\,.
$$
Thus, from \eqref{2.5}$_4$ and \eqref{2.12}--\eqref{2.14} it  follows that
$$\ba{rl}\medskip
\|\tilde{\bfN}(\bfu^{(0)}+\bfu^{(1)})\|&\!\!\!\!\le c_1\,\left[|\omega|\|\bfv\|_2+\|{\rm P}(\bfv\cdot\nabla\bfv)\|_2+|\omega|(|\gamma_1|+|\gamma_2|)\right]\\ \medskip
&\!\!\!\!\le c_2\,\left[(\|\bfu^{(1)}\|+\|\bfu^{(0)}\|)\|\bfu^{(1)}\| +\|\bfA_0^\alpha\bfv\|_2^2\right]\\
&\!\!\!\!\le c_3\,\left[(\|\bfu^{(1)}\|+\|\bfu^{(0)}\|)\|\bfu^{(1)}\| +\|\bfu^{(1)}\|_\alpha^2\right]\,,
\ea
$$
from which, using the boundedness of $\bfI^{-1}$, the validity of \eqref{H5} follows with $\kappa_1=\kappa_2=1$, and $\kappa_3=2$.\QED
\smallskip\par
As we know from the general approach, the stability property of the solution $\bfu\equiv\0$ to \eqref{1} requires that the eigenvalues in the spectral set $\sigma(\bfL_1)$ have all positive real part. 
The following \lemmref{2.1} provides the necessary and sufficient conditions for this to hold. However, its proof requires a simple but important preliminary result that we prove first.
\setcounter{lemm}{0}
\Bl For any ${\sf v}\in L^2(\mathscr C)$, we have
\be
\frac{{\sf C}_{\mathscr B}}{\sf C}\,\rho\,\|{\sf v}\|_2^2\le E:=\rho\|{\sf v}\|_2^2-{\sf C}\,a^2\le \rho\,\|{\sf v}\|_2^2\,.
\eeq{2.21}
\EL{2.0}
{\em Proof.} From \eqref{2.2}, we deduce
$$
E=\rho\,\|{\sf v}\|_2^2-\frac{\rho^2}{\sf C}\left(\int_{\mathscr C} (\bfe_3\times\bfx)\cdot{\sf v}\right)^2\,.
$$
Therefore, by the Schwarz inequality and \eqref{df}, we infer
$$
E\ge\rho\,\|{\sf v}\|_2^2-\frac\rho {\sf C}\left(\int_{\mathscr C} \rho (\bfe_3\times\bfx)^2\right)\|{\sf v}\|^2_2\ge \left(1-\frac{{\sf C}_{\mathscr L}}{{\sf C}}\right)\,\rho\,\|{\sf v}\|_2^2=\frac{{\sf C}_{\mathscr B}}{{\sf C}}\,\rho\,\|{\sf v}\|_2^2\,. 
$$
\QED
\par
We are now in a position to prove the following.
\Bl If  in  \eqref{2.4_1} $\xi=1$ then $\Re[\sigma(\bfL_1)]>0$, whereas if $\xi=-1$ then $\Re[\sigma(\bfL_1)]\cap (-\infty,0)\neq\emptyset$.
\EL{2.1}
{\em Proof.} To show the result it is sufficient to prove that all solutions to the equation
\be
\ode{\bfu}{t}+\bfL\bfu=\0\,,\ \ \bfu(0)\in H
\eeq{2.16}
are bounded, if $\xi=1$, whereas there exists at least one unbounded solution if $\xi=-1$. Now, \eqref{2.16} is equivalent to the following system of equations
\be\ba{ll}\medskip
\rho(\bfv_t+\dot{\omega}\,\bfe_3\times\bfx)-\mu\Delta\bfv=\nabla p\,,\ \ \ \nabla\cdot\bfv=0\,,\ \ \bfv|_{\partial\mathscr C}=\0\,,\\ \medskip
{\sf C}(\dot{\omega}-\dot{a})=\beta^2\gamma_2\,, \\ \medskip
\dot{\bfgamma}+\xi\,\omega\,\bfe_2=\0\,,\\ \left(\bfv(\cdot,0),\omega(0),\bfgamma(0)\right)\in L^2_\sigma(\mathscr C)\times\real\times\real^2\,.\ea
\eeq{2.17}
By dot-multiplying the first equation by $\bfv$, integrating by parts over $\mathcal C$ and employing \eqref{2.17}$_{2,3,4}$, we get
\be
\half\ode{}t\left[\rho\|\bfv\|_2^2-{\sf C}\,a^2\right]+\mu\|\nabla\bfv\|_2^2=\beta^2\gamma_2a\,.
\eeq{2.18}
Likewise, by multiplying both sides of \eqref{2.17}$_4$ by $\omega-a$ and using \eqref{2.17}$_5$, we deduce
\be
\half \ode{}t\left[{\sf C}\,(\omega-a)^2+\frac{\beta^2}\xi\,\gamma_2^2\right]=-\beta^2\gamma_2\,a\,.
\eeq{2.19}
As a result, summing \eqref{2.18} and \eqref{2.19} side by side we infer 
\be
\half\ode{}t\left[\rho\|\bfv\|_2^2-{\sf C}\,a^2+{\sf C}\,(\omega-a)^2+\frac{\beta^2}\xi\,\gamma_2^2\right]+\mu\|\nabla\bfv\|_2^2=0\,.
\eeq{2.20}
Consequently, if $\xi=1$, from \lemmref{2.0}, \eqref{2.20}, and \eqref{2.17}$_5$  we readily conclude that there is a constant $M>0$ depending on the initial data, such that
\be
\|\bfv(t)\|_2+|\omega(t)|+|\bfgamma(t)|\le M\,
\ \ \mbox{all $t\ge 0$}\,, \eeq{2.22}
which implies $\Re[\sigma(\bfL_1)]>0$. Suppose, next, $\xi=-1$ and, by contradiction, $\Re[\sigma(\bfL_1)]>0$. This means that any given solution to \eqref{2.16} must obey \eqref{2.22}. As a consequence, on the one hand,  from \eqref{2.18}, Schwarz and Poincar\'e inequalities and \eqref{2.21}  we get
\be
\ode{E}t+c_1\,E\le c_2\|\bfv\|_2\,.
\eeq{2.23}
On the other hand,
from \eqref{2.20} and again  Poincar\'e inequality, we infer
$$
\int_0^\infty \|\bfv(t)\|_2^2dt<\infty\,,
$$
so that the differential inequality \eqref{2.23} combined with \eqref{2.21} furnishes  (see \cite[Lemma 2.1]{GMZ})
$$
\lim_{t\to\infty}\|\bfv(t)\|_2=0\,.
$$
From the latter and \eqref{2.22} we easily deduce that the $\omega$-limit of the dynamical system generated by by \eqref{2.16} must be connected, compact and invariant, and therefore, in particular, that $\bfv\equiv\0$ there.
Using this property in \eqref{2.17}, we show also $\omega\equiv\gamma_2\equiv 0$. Thus, integrating  \eqref{2.20} from 0 to $t>0$, and then letting $t\to\infty$ we prove   
$$
2\mu\int_0^\infty\|\nabla\bfv(t)\|_2^2dt=\rho\|\bfv(0)\|_2^2-{\sf C}\,a^2(0)+{\sf C}\,(\omega(0)-a(0))^2+\frac{\beta^2}\xi\,\gamma_2^2(0)\,.
$$
However, since $\xi=-1$, this relation cannot be true for {\em any} initial data (it's enough to choose $\bfv(0)\equiv\0$ and $\omega^2(0)<(\beta^2/{\sf C})\gamma_2^2(0)$), thus showing a contradiction.
\QED
\par
Combining all we have shown so far in this section with \theoref{1}, we thus deduce the following stability results.
\Bt The steady-state solution ${\sf s}_0^+$ in \eqref{2.3}, representing the equilibrium configuration where the center of mass $G$ of $\mathscr S$ is in its lower position, is asymptotically, exponentially stable. Precisely,  the following properties hold. 
\begin{itemize}
  \item [{\rm (a)}] There is $\rho_0>0$ such that if, for some $\alpha\in [3/4,1)$, 
$$  
\|\bfA_0^\alpha\bfv(0)\|_2+|\omega(0)|+|\bfgamma(0)|<\rho_0\,,
$$
then there exists a corresponding  unique, global solution $(\bfv,\omega,\bfgamma)$ to \eqref{2.4} with $\xi=1$, such that,  
for all $T>0$,
$$\ba{ll}\medskip
\bfv\in C((0,T];D(\bfA_0))\cap C^1((0,T];L^2_\sigma(\mathcal C))\,,\ \ \bfA_0^\alpha\bfv\in C([0,T];L^2_\sigma(\mathcal C))\,,\\
\omega\in C([0,T];\real)\cap C^1((0,T];\real)\,;\ \  \bfgamma\in C^1([0,T];\real^2)\cap C^2((0,T];\real^2)\,;
\ea$$  
\item[{\rm (b)}] For any $\varepsilon>0$ there is $\delta>0$ such that
$$
\|\bfA_0^\alpha\bfv(0)\|_2+|\omega(0)|+|\bfgamma(0)|<\delta\ \ \Longrightarrow\ \ \sup_{t\ge 0}\left(\|\bfA_0^\alpha\bfv(t)\|_2+|\omega(t)|+|\bfgamma(t)|\right)<\varepsilon\,;
$$ 
\item[{\rm (c)}] There are $\eta,c,\kappa>0$ such that \end{itemize}
$$\ba{cc}\ms  
\|\bfA_0^\alpha\bfv(0)\|_2+|\omega(0)|+|\bfgamma(0)|<\eta\ \ \Longrightarrow\\ 
 \|\bfA_0^\alpha\bfv(t)\|_2+|\omega(t)|+|\bfgamma(t)|\le c\,\left(\|\bfA_0^\alpha\bfv(0)\|_2+|\omega(0)|+|\gamma_2(0)|\right)\,{\rm 
e}^{-\kappa\,t}\,,\ \mbox{all $t>0$}\,.
\ea$$\smallskip\par\noindent
Finally, the steady-state solution ${\sf s}_0^-$ in \eqref{2.3}, representing the equilibrium configuration where the center of mass $G$ of $\mathscr S$ is in its higher position, is unstable.
\ET{2}
{\em Proof.} In view of what we have already shown, the only thing that remains to prove is  $\bfgamma(t)\to \0$ as $t\to \infty$. To this end, we observe that, by \theoref{1}, $\bfgamma(t)\to\sigma\bfe_1$, for some $\sigma\in\real$, which,  by \eqref{fan} and property (b) above, must satisfy
$$  
\sigma^2+2\sigma=0\,,\ \ |\sigma|<\varepsilon\,.
$$
In turn, by choosing $\varepsilon$ (namely, $\delta$) appropriately, the latter implies $\sigma=0$, thus completing the proof.
\QED
\Br The {\em simple} stability of ${\sf s}_0^+$ as well as the instability  of ${\sf s}_0^-$, in the norm  $\|\,\|$ was  established in \cite[Theorem 1.5.2]{GaMa}, by a different method based on the study of the local dynamics of $\mathscr S$ near ${\sf s}_0^\pm$.
\ER{2}
\renewcommand{\theequation}{3.\arabic{equation}}\setcounter{equation}{0}\setcounter{section}{3}
\section*{\large 3. Asymptotic Behavior of the Motion of a Pendulum with a Liquid-Filled Cavity for Large Initial Data}
Another way of stating the stability part in \theoref{2} is to say  that all solutions to \eqref{2.1}--\eqref{2.2}  with ``sufficiently smooth" initial data that are  ``sufficiently close" to the equilibrium configuration ${\sf s}_0^+$ must remain ``close" to ${\sf s}_0^+$ and eventually converge to it at an exponential rate.  Objective of this section is to show that, in fact, the same conclusion holds in the more general  class of {\em weak} solutions to \eqref{2.1}--\eqref{2.2} and for data that not only are less regular, but also not necessarily ``close" to the stable equilibrium configuration ${\sf s}_0^+$. 
This result is achieved by suitably combining the findings of \cite{GaMa} with those of \theoref{2}. 
\setcounter{defn}{0}
\vspace*{.1mm}\par
We begin to recall the definition of weak solution \cite{GaMa}.
\Bd The triple $(\bfv,\omega,\bfchi)$ is a {\em weak solution} to \eqref{2.1} if it meets the following requirements:
\begin{itemize}
\item [{\rm (a)}]$\bfv\in C_w([0,\infty); L^2_\sigma(\mathscr C))\cap 
L^\infty(0,\infty;L^2_\sigma(\mathscr C))\cap L^2(0,\infty;W^{1,2}_0(\mathscr C))$\,;
\item [{\rm (b)}]$\omega\in C^0([0,\infty))\cap L^\infty(0,\infty)\,,\ \bfchi\in C^{1}([0,\infty);{\sf S}^1)\,;\,$\footnote{As customary, ${\sf S}^1$ denotes the unit sphere in $\real^2$.}
\item [{\rm (c)}] Strong Energy Inequality: namely, for all $t\ge s$ and a.a. $s\ge 0$ including $s=0$\,,
\be
\mathcal E(t)+\mathcal U(t)+\mu\Int{s}{t}\|\nabla \bfv(\tau)\|_2^2\,d\tau\le \mathcal E(s)+\mathcal U(s)
\eeq{sei}
where
$$
\mathcal E:=\half\big[\rho\,\|\bfv\|_2^2-{\sf C}\,a^2+{\sf C}\,(\omega-a\bigr)^2\bigr] \ \ \ \mbox{(kinetic energy)}
$$
and
$$
\mathcal U:=-{\sf C}\beta^2\chi_1\ \ \ \mbox{(potential energy)}
$$
\item[{\rm (d)}] $(\bfv,\omega,\bfchi)$ satisfies \eqref{2.1}$_{1,2,4,5}$ in the sense of distributions and \eqref{2.1}$_3$ in the trace sense.\end{itemize}\label{def}
\EED{1}
\par
The proof of the following important result is found in \cite[Proposition 1.3.6, Theorem 1.4.4]{GaMa}.
\Bp For any given initial data
\be
(\bfv_0,\omega_0,\bfchi_0)\in L^2_\sigma(\mathscr C)\times\real\times{\sf S}^1\,,
\eeq{3.1}
there exists at least one corresponding weak solution $(\bfv,\omega,\bfchi)$ and a time $t_0$ (depending on the solution) such that, setting $I_{t_0,T}=(t_0,t_0+T)$\,,
\be\ba{ll}\medskip 
\bfv\in C^0(\overline{I_{t_0,T}}; W_0^{1,2}(\mathscr C))\cap L^\infty(t_0,\infty;W_0^{1,2}(\mathscr C))\cap L^2(I_{t_0,T}; W^{2,2}(\mathscr C))\,,\\
\bfv_t\in L^2(I_{t_0,T}; H(\mathscr C)),\ \ \omega\in W^{1,\infty}(I_{t_0,T})\,,\ \ \bfchi\in W^{2,\infty}(I_{t_0,T};{\sf S}^1)\,,
\ea
\eeq{3.2}
for all $T>0$. Moreover, there is $p\in L^2(I_{t_0,T};W^{1,2}(\mathscr C))$, all $T>0$, such that $(\bfv,p,\omega,\bfgamma)$ satisfies \eqref{2.1}$_{1,2}$ a.e. in $\mathscr C\times (t_0,\infty)$. In addition, the following asymptotic properties hold:
$$
\lim_{t\to\infty}\left(\|\bfv(t)\|_{2,2}+\|\bfv_t(t)\|_2+|\omega(t)|\right)=0\,.
$$
Finally, for all initial data such that
\be
\rho\,\|\bfv_0\|_2^2+{\sf C}\,(\omega_0^2-a(0))^2<2\,{\sf C}\,\beta^2\,(1+ \chi_{1,0})\,
\eeq{3.3}
we have also
$$
\lim_{t\to\infty}|\bfchi(t)-\bfe_1|=0\,,
$$
that is, the coupled system pendulum-liquid goes to the equilibrium configuration with the center of mass in its lower position.
\EP{1}\par
From \propref{1} and \theoref{2} we are now able to prove the  main result of this section that establishes the rate of decay to equilibrium.\setcounter{theo}{0}
\Bt Let the initial data \eqref{3.1} satisfy condition \eqref{3.3}. Then, for any corresponding weak solution $(\bfv,\omega,\bfchi$), there are $t_0,C_1$, possibly depending on the solution, and $C_2>0$ such that
$$
\|\bfv(t)\|_{2,2}+\|\bfv_t(t)\|_2+|\omega(t)|+|\dot{\omega}(t)|+|\bfchi(t)-\bfe_1|\le C_1\,{\rm e}^{-C_2\,t}\,,\ \ \mbox{for all $t\ge t_0$}\,.
$$
\ET{3}
{\em Proof.} Recalling that \cite[Section IV.6]{Ga} 
$$
\|\bfv\|_{2,2}\le C_3\|\bfA_0\bfv\|_2\le C_4\|\bfv\|_{2,2}\,,
$$
and $\bfchi=\bfgamma+\bfe_1$, 
by \propref{1} and \eqref{2i} (with $\bfA\equiv\bfA_0$,\ $\|\,\|\equiv\|\,\|_2$) we get that there exists $t_0>0$ such that
$$
\|\bfA_0^\alpha\bfv(t_0)\|_2+|\omega(t_0)|+|\bfchi(t_0)-\bfe_1|<\eta
$$
with  $\eta$ as in \theoref{2}(c), and $\alpha\in [0,1]$. Therefore, again by \theoref{2}(c) and \eqref{2i} we deduce
\be
\|\bfA_0^\alpha\bfv(t)\|_{2}+|\omega(t)|+|\bfchi(t)-\bfe_1|\le C\,{\rm e}^{-\kappa\,t}\,,\ \ \mbox{for all $\alpha\in [0,1)$ and $t\ge t_0$}\,.
\eeq{3.5}
Now, from \cite[eq. (1.44)]{GaMa} we know that, for all $t\ge t_0$, the following differential inequality holds
\be
\ode{}tE_1+c_1\|\nabla\bfv_t\|_2\le c_2\left[\omega^2+\|\bfv\|_2+(\|\bfv\|_2+\|\nabla\bfv\|_2^8)\|\bfv_t\|_2^2\right]\,,
\eeq{3.6}
where
\be
\frac{{\sf C}_{\mathscr B}}{{\sf C}}\rho\|\bfv_t\|_2^2\le E_1:=\rho\,\|\bfv_t\|_2^2-\frac{\rho^2}{{\sf C}}\left(\int_{\mathscr C}(\bfe_3\times\bfx)\cdot\bfv_t\right)^2\le\rho \|\bfv_t\|^2\,;
\eeq{3.7}
see \lemmref{2.0}. Therefore, employing the Poincar\`e inequality $\|\nabla\bfv_t\|_2\ge c\,\|\bfv_t\|_2$, and taking into account \eqref{3.5} and \eqref{3.7}, from \eqref{3.6} we deduce
$$
\ode{}tE_1+c_3\,E_1\le c_4{\rm e}^{-c_5t}(1+E_1)\,. 
$$
By a direct application of a Gronwall-like lemma to the latter inequality, and also with the help of \eqref{3.7}, we infer
\be
\|\bfv_t(t)\|_2\le c_6\,{\rm e}^{-c_7t}\,,\ \ t\ge t_0\,.
\eeq{3.8}
Plugging this information back in \eqref{2.4}$_4$ and using \eqref{3.5} entails
\be
|\dot{\omega}(t)|_2\le c_8\,{\rm e}^{-c_9t}\,,\ \ t\ge t_0\,.
\eeq{3.9}
Next, employing \eqref{2.1}$_1$ and \eqref{2.1}$_4$, one can show the following estimate \cite[eq. (1.46)]{GaMa}
\be
\|\bfv(t)\|_{2,2}\le c_{10}\left(\|\nabla\bfv(t)\|_{2}^3+\|\bfv_t(t)\|_2+|\omega(t)|\,\|\bfv(t)\|_2+|\chi_2(t)|\right)\,,\ \ t\ge t_0\,. 
\eeq{3.10}
The desired result is then a consequence of \eqref{3.5}, and \eqref{3.8}--\eqref{3.10}.\QED
{\bf Acknowledgment.} The work of G.P.~Galdi is partially supported by NSF grant DMS-1614011

\ed

\